\documentclass[11pt]{article}
\usepackage[usenames]{color}
\usepackage{amssymb,psfig,epsfig,latexsym,graphicx}
\usepackage{amsmath}
\usepackage{ifpdf}

\usepackage[english]{babel}
\usepackage{setspace}
\usepackage{dsfont}

\definecolor{webgreen}{rgb}{0,.5,0}
\definecolor{webbrown}{rgb}{.6,0,0}

\newcommand{\eqn}[1]{(\ref{#1})}

\newcommand{\beql}[1]{\begin{equation}\label{#1}}
\newcommand{\eeq}{\end{equation}}
\newcommand{\arr}{\leftrightarrow}

\DeclareMathOperator{\divides}{divides}

\newcommand{\sU}{{\mathcal U}}
\newcommand{\sE}{{\mathcal E}}
\newcommand{\sV}{{\mathcal F}}
\newcommand{\sZ}{{\mathcal Z}}
\newcommand{\sN}{{\mathcal N}}

\newsavebox{\gplussave}
\sbox{\gplussave}{\raisebox{0ex}{\includegraphics[width=8pt]{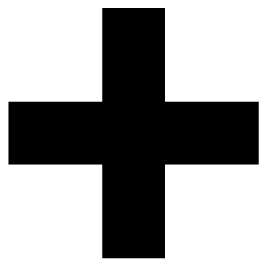}}}
\newcommand{\gplus}{\usebox{\gplussave}}

\newsavebox{\gminussave}
\sbox{\gminussave}{\raisebox{0ex}{\includegraphics[width=8pt]{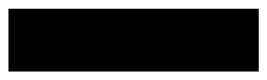}}}
\newcommand{\gminus}{\usebox{\gminussave}}

\newsavebox{\gtimessave}
\sbox{\gtimessave}{\raisebox{-.1ex}{\includegraphics[width=9pt]{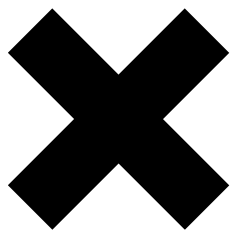}}}
\newcommand{\gtimes}{\usebox{\gtimessave}}

\newcommand{\CA}{\,\gplus\,}
\newcommand{\CS}{\,\gminus\,}
\newcommand{\CM}{\,\gtimes\,}

\DeclareMathOperator{\Rim}{Rim}

\makeatletter
\def\@sect#1#2#3#4#5#6[#7]#8{\ifnum #2>\c@secnumdepth
 \def\@svsec{}\else
 \refstepcounter{#1}\edef\@svsec{\csname the#1\endcsname.\hskip .75em }\fi
 \@tempskipa #5\relax
 \ifdim \@tempskipa>\z@
 \begingroup #6\relax
 \@hangfrom{\hskip #3\relax\@svsec}{\interlinepenalty \@M #8\par}%
 \endgroup
 \csname #1mark\endcsname{#7}\addcontentsline
 {toc}{#1}{\ifnum #2>\c@secnumdepth \else
 \protect\numberline{\csname the#1\endcsname}\fi
 #7}\else
 \def\@svsechd{#6\hskip #3\@svsec #8\csname #1mark\endcsname
 {#7}\addcontentsline
 {toc}{#1}{\ifnum #2>\c@secnumdepth \else
 \protect\numberline{\csname the#1\endcsname}\fi
 #7}}\fi
 \@xsect{#5}}
\def\@begintheorem#1#2{\it \trivlist \item[\hskip \labelsep{\bf #1\ #2.}]}

\def\section{\@startsection {section}{1}{\z@}{-3.5ex plus -1ex minus
 -.2ex}{2.3ex plus .2ex}{\normalsize\bf}}
\makeatother
\makeatletter
\def\subsection{\@startsection {subsection}{1}{\z@}{-3.5ex plus -1ex minus
 -.2ex}{2.3ex plus .2ex}{\normalsize\bf}}

\makeatother

\begin{document}

\begin{center}
{\Large\bf Carryless Arithmetic Mod 10 } \\
\vspace*{+.2in}
David~Applegate, Marc~LeBrun and N.~J.~A.~Sloane
\end{center}

\vspace*{+.5in}

\section*{Nim}\label{Sec0}

Forms of Nim have been played since antiquity and a complete theory was 
published as early as 1902 (see \cite{Bout02}).  
Martin Gardner described the game in one of his earliest columns \cite{MG59}
and returned to it many times over the years 
(\cite{MG77}--\cite{MG01}).

Central to the analysis of Nim is Nim-addition. 
The Nim-sum is calculated by writing the terms in base $2$ 
and adding the columns mod $2$, with {\em no carries}.
A Nim position is a winning position if and only if the Nim-sum
of the sizes of the heaps is zero \cite{WW}, \cite{MG59}.

Is there is a generalization of Nim in which
the analysis uses the base-$b$ representations of
the sizes of the heaps, for $b>2$,
in which a position is a win if and only if the mod-$b$
sums of the columns is identically zero?
One such game, $\Rim_b$ (an abbreviation of Restricted-Nim) exists, 
although it is complicated and not well known. 
It was introduced in an unpublished paper \cite{Flan80}
in 1980 and is hinted at in \cite{Ferg98}.
Despite his interest in Nim, Martin Gardner never mentions $\Rim_b$, nor
does it appear in {\em Winning Ways} \cite{WW},
which extensively analyzes Nim variants.

In the present paper we focus on $b=10$,
and consider, not $\Rim_{10}$ itself, but the arithmetic
that arises if calculations, addition and multiplication,
are performed mod $10$, with {\em no carries}. 
Along the way we encounter several new and interesting
number sequences, which would have appealed to
Martin Gardner, always a fan of integer sequences.

\vspace*{+.1in}

\section*{The Carryless Islands}\label{Sec1}

The fabled carefree residents of the 
Carryless Islands in the remote South Pacific
have very few possessions, which is just as well, since their
arithmetic is ill-suited to accurate bookkeeping.
When they add or multiply numbers, they follow rules similar
to ours, except that there are {\em no carries} into
other digit positions. Sociologists
explain this by noting that the Carryless Islands
were originally penal colonies, and, as penal institutions are
generally known to have excellent dental care, the
islanders were, happily, generally free of carries.
We will use $\CA$ and $\CM$ for their operations,\footnote{We
prefer not to use outlandish symbols such as $\clubsuit$
and $\spadesuit$, since $\CA$ and $\CM$
are perfectly reasonable operations,
although to our eyes they have rather strange properties.
As Marcia Ascher remarks, writing about
mathematics in indigenous cultures,
``in many cases these cultures and their ideas were
unknown beyond their own boundaries, or misunderstood when first encountered
by outsiders'' \cite{Asc}.}
and $+$ and $\times$ for the standard operations used
by the rest of the world.
Addition and multiplication of single-digit numbers
are performed by ``reduction mod $10$.''
Carry digits are simply ignored, so 
$9 \CA 4=3$,
$5 \CA 5=0$,
$9 \CM 4=6$,
$5 \CM 4=0$,
and so on.
Adding or multiplying 
larger numbers also follows the familiar procedures, but again 
with the proviso that there are no carries.
For example, adding $785$ and $376$ produces 51,
and the product of $643$ and $59$ is $417$ (see Figure 1).

  \begin{minipage}[t]{2.5in}
\begin{center}
$$
\begin{tabular}{ p{.001in}  p{.001in}  p{.001in}  p{.001in} }
       & 7 & 8 & 5 \\
$\CA$  & 3 & 7 & 6 \\
\hline
       & 0 & 5 & 1 \\
\hline
%\caption{Fig. 1(a)}\label{Fig1a}
\end{tabular}
$$
Fig. 1(a) Carryless addition.
\end{center}
  \end{minipage}
  \begin{minipage}[t]{2.5in}
\begin{center}
$$
\begin{tabular}{ p{.001in} p{.001in}  p{.001in}  p{.001in} }
       & 6 & 4 & 3 \\
$\CM$  &   & 5 & 9 \\
\hline
       & 4 & 6 & 7 \\
   0 & 0 & 5 &   \\
\hline
   0 & 4 & 1 & 7 \\
\hline
%\caption{Fig. 1(b)}\label{Fig1b}
\end{tabular}
$$
Fig. 1(b) Carryless multiplication. 
\end{center}
  \end{minipage}
\vspace*{+.2in}

What does elementary number theory look like on these islands?
Let's start with the carryless squares $n \CM n$. For $n=0,1,2,3$
we get $0$, $1$, $4$, $9$, 
Then for $n > 3$ we have
$4 \CM 4 = 6$,
$5 \CM 5 = 5$,
$6 \CM 6 = 6$,
$7 \CM 7 = 9$,
$8 \CM 8 = 4$,
$9 \CM 9 = 1$,
$10 \CM 10 = 100, \ldots$,
giving the sequence
$$
0,\, 1,\, 4,\, 9,\, 6,\, 5,\, 6,\, 9,\, 4,\, 1,\, 100,\, 121,\, 144,\, 169,\, 186,\, 105,\, 126,\, 149,\, 164,\, \ldots\, .
$$ 
It turns out that this is entry A059729 in the OEIS \cite{OEIS},
contributed by Henry Bottomley on February 20, 2001,
although without any reference to earlier work on these numbers.
Bottomley also contributed sequence A059692, giving the 
carryless multiplication table, and several other sequences
related to carryless products.
Likewise the sequence of values of $n \CA n$, 
$$
0,\, 2,\, 4,\, 6,\, 8,\, 0,\, 2,\, 4,\, 6,\, 8,\, 20,\, 22,\, 24,\, 26,\, 28,\, 20,\, 22,\, 24,\, 26,\, 28,\, 40,\, 42,\,  \ldots\,,
$$
is entry A004520, submitted to the OEIS by one of the present authors
around 1996, again without references.
(If these numbers are sorted and duplicates removed, we get 
the carryless ``evenish'' numbers, that is, numbers
all of whose digits are even, A014263.)
Carryless arithmetic must surely have been studied before now, but
the absence of references in \cite{OEIS} suggests that it is
not mentioned in any of the standard texts on number theory.

\section*{The carryless primes}\label{Sec3}

If we require that a prime $\pi$ is a number
whose only factorization is $1$ times itself, we are out
of luck, since every carryless number is divisible by $9$,
and there would be no primes at all. (For $9 \CM 1 = 9$, 
$9 \CM 2 = 8$, $9 \CM 3 = 7, \ldots$, $9 \CM 9 = 1$. 
So if we construct a number $\rho$ by replacing all the $1$'s in $\pi$ by $9$'s,
all the $2$'s by $8$'s, $\ldots$ then $\pi= 9 \CM \rho$,
and $\pi$ would not be a prime.)

There {\em are} primes, when defined in the right way.
Since $1 \CM 1 = 1$, $3 \CM 7 = 1$ and $9 \CM 9=1$,
all of $1,3,7$ and $9$ divide $1$ and so divide any number.
We call $1, 3, 7$ and $9$ {\em units},
the usual name for integers that divide $1$.
Units should not be counted as factors when considering if a number is prime
(just as factors of $-1$ are ignored in ordinary
arithmetic: $7 = (-1) \times (-7)$ doesn't count as a
factorization when considering if $7$ is a prime).

So we define a carryless prime to be a non-unit $\pi$
whose only factorizations are of the form $\pi = u \CM \rho$
where $u$ is a unit.
Computer experiments suggest that the first few primes are
\beql{Eq0}
21,\, 23,\, 25,\, 27,\, 29,\, 41,\, 43,\, 45,\, 47,\, 49,\, 51,\, 52,\, 53,\, 54,\, 56,\, 57,\, 58,\, 59,\, 61,\, 63,\, \ldots \,,
\eeq
but there are surprising omissions in this list, resulting
from some strange factorizations: $2 = 2 \CM 51$, $10 = 56 \CM 65$, $11 = 51 \CM 61$.
It is hard to be sure at this stage that the above list
is correct, since there exist factorizations
where one of the numbers is much larger than the number being factored,
such as $2 = 4 \CM 5005505553$.
One property that makes carryless arithmetic interesting is the 
presence of {\em zero-divisors}: the product of two numbers can be zero without
either of them being zero: $2 \CM 5 = 0$, $628 \CM 55 = 0$.
Perhaps $21$ is the product of two really huge numbers?
Nonetheless, the list {\em is} correct, as we will see 
(it is now sequence A169887 in \cite{OEIS}).

\section*{Algebra to the rescue}\label{Sec4}

The secret to understanding carryless 
arithmetic is to introduce a little algebra.
Let $R_{10}$ denote the ring of integers mod $10$, and $R_{10}[X]$
the ring of polynomials in X with coefficients in $R_{10}$.
Then we can represent carryless numbers by elements of $R_{10}[X]$:
$21$ corresponds to $2X+1$, $109$ to $X^2+9$, and so on.
Carryless addition and multiplication are
simply addition and multiplication in $R_{10}[X]$:
our first example,
$$
785 \CA 376 = 51 \,.
$$
corresponds to 
$$
(7X^2+8X+5) + (3X^2+7X+6) = 5X+1\,,
$$
where the polynomials are added or multiplied 
in the usual way, and the coefficients then reduced mod $10$.
Conversely, any element of $R_{10}[X]$ represents a 
unique carryless number (just set $X=10$ in the polynomial). 
In fact arithmetic in $R_{10}[X]$
is clearly exactly the same
as the arithmetic
of carryless numbers.
This could be used
as a formal definition of carryless arithmetic mod $10$.
It also shows that this
arithmetic is commutative, associative and distributive.

Since $R_{10}[X]$ is a ring, we can not only add and multiply,
we can also subtract, something the Carryless Islanders never
considered.
The negatives of the elements of $R_{10}$ are
$-1 = 9$,
$-2 = 8, \ldots,$
$-9 = 1$,
and similarly for the elements of $R_{10}[X]$.
So the negative of a carryless number is
its ``10's complement,''
obtained by replacing each nonzero digit $d$ by $10-d$,
for example $\CS 702 = 308$.
To subtract $A$ from $B$, we add $\CS A$ to $B$:
$650 \CS 702 = 650 \CA 308 = 958$.
This is equivalent to doing elementary school subtraction
where we can ``borrow'' but don't have to pay back!

The units in $R_{10}[X]$, that is, the elements that divide $1$,
are the constants $1,3,7,9$, 
and the carryless primes that we defined
are the {\em irreducible} elements in $R_{10}[X]$, that is, 
non-units $f_{10}(X) \in R_{10}[X]$
whose only factorizations are of the form
$f_{10}(X)=u\, g_{10}(X)$, where $u$ is a unit and $g_{10}(X) \in R_{10}[X]$.
The units can also be written as $1, \CS 1, 3$ and $\CS 3$,
which more closely relates them to the units $1$ and
$-1$ in ordinary arithmetic ($3$ and $\CS 3$ act
in some ways like the imaginary units $i$ and $-i$, squaring to $-1$, for example).

The key to further progress is to notice that $R_{10}$
is the direct sum of the ring
$R_2$ of integers mod $2$ and the ring $R_5$ 
of integers mod $5$. Given $r_{10} \in R_{10}$, we read it mod $2$ 
and mod $5$ to obtain a pair $[r_2, r_5]$
with $r_2 \in R_2$, $r_5 \in R_5$.
The elements $0,1,\ldots,9 \in R_{10}$ (or equivalently the
carryless digits $0,1,\ldots,9$) 
and their corresponding pairs $[r_2, r_5]$ are given
by the following table. The Chinese Remainder Theorem guarantees
that this is a one-to-one correspondence.
\beql{Eq1}
\begin{array}{cccccccccc}
0 & 1 & 2 & 3 & 4 & 5 & 6 & 7 & 8 & 9   \\
$[0,0]$ & $[1,1]$ & $[0,2]$ & $[1,3]$ & $[0,4]$ & $[1,0]$ & $[0,1]$ & $[1,2]$ & $[0,3]$ & $[1,4]$ 
\end{array}
\eeq
As a check, we note that $\{1\}$ is the (singleton) set of units in $R_2$,
while $\{1,2,3,4\}$ is the set of units in $R_5$, so the pairs 
$[1,1]$, $[1,2]$, $[1,3]$
and $[1,4]$ correspondingly produce the units $1$, $7$, $3$ and $9$ of $R_{10}$.

Similarly, polynomials $f_{10}(X) \in R_{10}[X]$ correspond
to pairs of polynomials $[f_2(X), f_5(X)]$,
obtained by reading $f_{10}(X)$ respectively mod $2$ and mod $5$.
Conversely, given any such pair of polynomials
$[f_2(X), f_5(X)]$, there is a unique $f_{10}(X) \in R_{10}[X]$
that corresponds to them, which can be found using \eqn{Eq1}.
We indicate this by writing $f_{10}(X) \arr [f_2(X),f_5(X)]$.
If also $g_{10}(X) \arr [g_2(X),g_5(X)]$, then
$f_{10}(X)+g_{10}(X) \arr [f_2(X)+g_2(X),f_5(X)+g_5(X)]$ and
$f_{10}(X)g_{10}(X) \arr [f_2(X)g_2(X),f_5(X)g_5(X)]$.

We are now in a position to answer many questions
about carryless arithmetic.

\section*{The carryless primes, again}\label{Sec5}

What are the irreducible elements $f_{10}(X) \in R_{10}[X]$?
If $f_{10}(X) \arr [f_2(X),f_5(X)]$ is irreducible then certainly 
$f_2$ and $f_5$ must be either units or irreducible,
for if $f_2=g_2 h_2 $ then we have the 
factorization $[f_2,f_5] = [g_2,f_5] [h_2,1]$.
Also $[f_2,f_5]=[f_2,1][1,f_5]$,
so one of $f_2$, $f_5$ must
be irreducible and the other must be a unit.
So the irreducible elements
in $R_{10}[X]$ are of the form $[f_2(X),u]$, where $f_2(X)$
is an irreducible polynomial mod $2$ of degree $\ge 1$
and $u \in \{1,2,3,4\}$,
together with elements of the form $[1,f_5(X)]$, where $f_5(X)$
is an irreducible polynomial mod $5$ of degree $\ge 1$.

The irreducible polynomials mod $2$ are $X$, $X+1$, $X^2+X+1, \ldots$,
and 
the irreducible polynomials mod $5$ are 
$uX$, $uX+v, \ldots$,
where $u, v \in \{1,2,3,4\}$
(see entries A058943, A058945 in \cite{OEIS}).
The first few irreducible elements in $R_{10}[X]$ are therefore
$[X,1]$, 
$[X,2]$, 
$[X,3]$, 
$[X,4]$, 
$[X+1,1]$, 
$[X+1,2], \ldots$, 
and
$[1,X]$, 
$[1,2X]$, 
$[1,3X]$, 
$[1,4X]$, 
$[1,X+1]$,
$[1,2X+1], \ldots$.
The corresponding carryless primes, according to \eqn{Eq1}, are
$56, 52, 58, 54, 51, 57, \ldots,$ and
$65, 25, 85, 45, 61, 21, \ldots$.
And so we can verify that the list in \eqn{Eq0} is correct.

We will call a number with at least two digits in which all digits
except the rightmost are even but the rightmost is odd
an {\em e-type number} (A143712), and a number 
with at least two digits in which all digits except the rightmost are $0$ or $5$
and the rightmost is neither $0$ nor $5$ an {\em f-type number} (A144162).
Similarly, we call the 
primes corresponding to the irreducible elements $[1,f_5(X)]$
{\em e-type primes}, and the primes
corresponding to the irreducible elements 
$[f_2(X),u]$ {\em f-type primes}.

We also see that our earlier concern about the primality of $21$
was groundless. It is impossible
for the length (in decimal digits) of a
nonzero carryless product to be less than
the length of both of the factors. This follows
from the fact that if $\ell(n)$ is
the number of
decimal digits in the number $n>0$ corresponding to a pair
$[f_2(X), f_5(X)]$, then
$\ell(n) = 1 + \max\{\deg f_2, \deg f_5\}$.
So if $mn > 0$, $\ell(mn) \ge \min\{\ell(m), \ell(n)\}$.

Also, since we know how many irreducible polynomials
mod $2$ and mod $5$ there are of given degree 
(see A001037, A001692 in \cite{OEIS}),
we can write down a formula for the number of
$k$-digit carryless primes, something that 
we cannot do for ordinary primes, namely
$$
\frac{4}{k-1} \sum_{d \divides k-1} \mu \left( \frac{k-1}{d} \right)
(2^d + 5^d) \,,
$$
for $k \ge 2$, where $\mu$ is the M\"obius function (A008683).
There are $28$ primes with two digits
(the twenty listed in \eqn{Eq0}, together
with $65, 67, 69, 81, 83, 85, 87, 89$),
 $44$ with three digits, $\ldots$ (A169962).
For large $k$ the number is about $4 \cdot 5^{k-1}/(k-1)$,
whereas the number of ordinary primes
with exactly $k$ digits is much larger, about $9\cdot 10^{k-1}/(k \log 10)$,
so carryless primes are much rarer than ordinary primes.

Incidentally, the {\em prime ideals}
in $R_{10}[X]$, as distinct from the irreducible elements,
all have a single generator, which is one of
$[0,1], [1,0], [1,1], [f_2(X),1], [1, f_5(X)]$,
where $f_2(X)$, $f_5(X)$
are irreducible (cf. \cite[Chap.~III,~Thm.~30]{ZS1}).

\section*{The carryless squares, again}\label{Sec6}

Squaring a mod $2$ polynomial is easy: $f_2(X)^2 = f_2(X^2)$.
So if $n$ corresponds to the pair
$[f_2(X), f_5(X)]$, $n^2$ corresponds to $[f_2(X^2), f_5(X)^2]
= [f_2(X^2),0] + [0, f_5(X)^2]$.
This gives a two-step recipe for producing all carryless squares.
First find (using \eqn{Eq1}) the
carryless number $m$ corresponding to $[0, f_5(X)^2]$,
where $f_5(X)$ is any polynomial mod $5$.
The effect of adding a nonzero $[f_2(X^2),0]$ changes
some subset of
the digits in positions $0, 2, 4, \ldots$ of $m$ by the addition of $5$ mod $10$.

For example, if $f_5(X) = X+2$, $f_5(X)^2=X^2+4X+4$,
and by \eqn{Eq1} $[0, f_5(X)^2]$ corresponds to the carryless square $m = 644$.
We now add $5$ mod $10$ to any subset
of the digits in positions $0, 2, 4, 6, \ldots$ of $m$
(considering $m$ extended by prefixing it with any number of zeros),
obtaining infinitely many squares
$644, 649, 144, 149, \ldots, 50644, 5050649, \ldots$.

This also leads to a formula for the number of $k$-digit
carryless squares. For even $k$ the number is $0$,
and for odd $k$ it is
$$
\frac{1}{2} \, 9 \cdot 10^{(k-1)/2} ~+~ 2^{(k-3)/2} 
$$
(zero is excluded from the count).
There are five squares of length $1$ (namely $1, 4, 5, 6$ and $9$),
$46$ of length $3$, $\ldots$
(see A059729, A169889, A169963).
For large odd $k$ there are about twice as many $k$-digit 
carryless squares as ordinary squares.

\section*{Divisors and factorizations}\label{Sec7}

What about the factorization of numbers into
the product of carryless primes?
Unfortunately, the existence of zero-divisors complicates matters,
and it turns out that there is no natural way to
define, for example, an analog of the usual
sum-of-divisors function $\sigma(n)$.
In our analysis we define several classes of carryless numbers:
\begin{list}{}{\setlength{\itemsep}{0.02in}}
\item
$\sU := \{1,3,7,9\}$, the units,
\item
$\sE = \{0,2,4,6, 8, 20,22,\ldots\}$, the ``evenish'' numbers, 
in which all digits are even (A014263),
\item
$\sV = \{0,5,50,55,\ldots\}$, the ``fiveish'' numbers, 
in which all digits are $0$ or $5$ (A169964),
\item
$\sZ := \sE \cup \sV = \{0,2,4,5,6,8,20,22,\ldots\}$, the zero-divisors 
(A169884),
\item
$\sN = \{1,3,7,9,10,11,12,13,\ldots\}$, the positive numbers not in $\sZ$ 
(A169968).
\end{list}

Suppose $d$ is a carryless divisor of $n$, that is,
there is a number $q$ such that $d \CM q = n$.
What can be said about the possible choices for $q$?
One can show---we omit the straightforward proofs---that
\begin{list}{\textbullet}{\setlength{\itemsep}{0.02in}}
\item
if $d \in \sN$ then there is a unique $q$,
\item
if $d \in \sE$ then $d \CM q' = n$ if and only if $q'=q ~ \CA ~ v$ 
for some $v \in \sV$,
\item
if $d \in \sV$ then $d \CM q' = n$ if and only if $q'=q ~ \CA ~ e$ 
for some $e \in \sE$.
\end{list}

The same distinctions are needed to describe factorizations
into primes.
\begin{list}{\textbullet}{\setlength{\itemsep}{0.02in}}
\item
If $n \in \sN$ then $n$ has a unique factorization
as a carryless product of primes, up to 
multiplication by units.
For example, we already saw $10 = 56 \CM 65$.
But we also have $10 = (3 \CM 56) \CM (7 \CM 65)
= 58 \CM 25 = (9 \CM 56) \CM (9 \CM 65) = 54 \CM 45 = 9 \CM 52 \CM 25$, etc.,
illustrating the nonuniqueness. Also
$11 = 51 \CM 61$; $101 = 21 \CM 29 \CM 51$,
$1234 = 23 \CM 23 \CM 23 \CM 51 \CM 51 \CM 52$.
It follows that any 
non-unit in $\sN$ can be written both as
$e \CM f$ and $e' \CA f'$, where $e$ and $e'$ are e-type numbers
and $f$ and $f'$ are f-type numbers.
For example, $12 = 81 \CM 52 = 61 \CA 51$.
\item
If $n \in \sE$ then $n$ has a unique factorization
as $2$ times a product of e-type primes, up to multiplication by units
(in this case, every f-type prime divides $n$).
For example, $20 = 2 \CM 65$,
$22 = 2 \CM 61$, $2468 = 2 \CM 69 \CM 69 \CM 69$.
\item
If $n \in \sV$ then $n$ has a unique factorization
as $5$ times a product of f-type primes, up to multiplication by units
(in this case, every e-type prime divides $n$).
For example, $50 = 5 \CM 52$, $505 = 5 \CM 51 \CM 51$.
\end{list}

Here are the analogous statements about divisors:
\begin{list}{\textbullet}{\setlength{\itemsep}{0.02in}}
\item
if $n \in \sN$, $n$ has only finitely many divisors.
If $ d$ divides $n$ and $u \in \sU$, then $d \CM u$
divides $n$. The divisors may
be grouped into equivalence classes $d \CM \sU$.
Since the sum of the elements of $\sU$ is zero, so is
the sum of the divisors of $n$.
\item
if $n \in \sE$, $d$ divides $n$, $u \in \sU$ and $v \in \sV$, 
then $d \CM u ~ \CA ~ v$ divides $n$. So $n$ has
infinitely many divisors, belonging to
equivalence classes $d \CM \sU ~ \CA ~ \sV$.
\item
if $n \in \sV$, $d$ divides $n$, $u \in \sU$ and $e \in \sE$, 
then $d \CM u ~ \CA ~ e$ divides $n$. So $n$ has
infinitely many divisors, belonging to
equivalence classes $d \CM \sU ~ \CA ~ \sE$.
\end{list}

Any attempt to define a sum-of-divisors function must
specify how to choose representatives from the 
equivalence classes. There seems to be no natural way to do this.
One possibility would be to choose the smallest
decimal number in each class, but this seems 
unsatisfactory (since it depends on the ordering
of decimal numbers, another concept the islanders
seem not to be familiar with).

\section*{Further number theory}\label{Sec8}

In summary, we can help the Carryless Islanders
by defining subtraction, prime numbers, and factorization
into primes.  But further concepts such as
the number of divisors, the sum of divisors
and perfect numbers seem to lie beyond these Islands.

However, many other carryless analogs are well-defined, including
including triangular numbers (A169890),
cubes (A169885), partitions (A169973), greatest common
divisors and least common multiples, and so on.
Some seem exotic, while other familiar sequences 
simply become periodic. For example, the analog
of the Fibonacci numbers
coincides with the sequence of Fibonacci numbers read mod $10$,
A003893, which becomes periodic with period $60$
(the periodicity of the Fibonacci numbers to any
modulus being a well-studied subject, see sequence A001175).
Similarly, the analogue of the powers of $2$ (A000689) 
becomes periodic with period $4$.
We might also generalize beyond simple squares, cubes, etc. and
investigate the properties of polynomials or power series based on carryless
operations--How do these factor?  What are their fixed points?--and so on.

Taking a different tack, carryless mod $10$ partitions are enumerated
in A169973, which may be derived as the coefficients of $z^n$ in the formal
expansion of the analog of the classic partition generating function
$ \prod_{k=1}^{\infty} (1+z^k)$,
wherein powers of $z$ are multiplied together by combining their exponents
with carryless mod $10$ addition instead of the ordinary sum.

\section*{Afterword}\label{Sec9}
There's a great deal yet to be explored in these Carryless Islands!
Watch for our next paper on another carryless arithmetic, 
in which operations on single digits are defined by
$a \oplus b=\max \{a, b\}$,
$a \otimes b=\min \{a, b\}$.
We call this ``dismal arithmetic.''

When the {\em Handbook of Integer Sequences} was published
39 years ago, Martin Gardner was kind enough to write in his 
{\em Mathematical Games} column of July 1974 that 
``every recreational mathematician should buy
a copy forthwith.'' That book contained 2372 sequences:
today its successor, the {\em On-Line Encyclopedia of Integer Sequences}
(or OEIS) \cite{OEIS}, contains nearly 200,000 sequences.
We were about to write to Martin about carryless arithmetic
when we heard the sad news of his death. 
This article, the first of a series on various kinds of carryless arithmetic,
is offered in his honor.

\small
{\bf Summary}
What might arithmetic look like on an island that eschews carry digits?
How would primes, squares and other number theoretical concepts play out on such
an island?

\vspace*{+.5in}

{\bf David Applegate} (\texttt{david@research.att.com}) received his
Ph.D. in Computer Science from Carnegie Mellon, and has been at AT\&T
Shannon Labs since 2000.  His research interests include combinatorial
optimization, the Traveling Salesman Problem, and other mathematical
diversions.

\vspace*{+.5in}

{\bf Marc LeBrun} was confirmed as an ardent amateur mathematician at age ten,
when Martin Gardner graciously answered his fan letter.  He is delighted to
find professional colleagues who share his interests in recreational
arithmetic, and is an enthusiastic contributor to the {\em On-Line Encyclopedia of
Integer Sequences}.

\vspace*{+.5in}

{\bf N. J. A. Sloane} (home page \texttt{www.research.att.com/$\sim$njas}) has been
at AT\&T Bell Labs and later AT\&T Shannon Labs since 1969,
and has written many books and articles on mathematics, engineering and statistics.
He also runs the {\em On-Line Encyclopedia of Integer Sequences}.

\vspace*{+.5in}

AMS 2010 Classification: Primary 06D05, 11A63

\end{document}